\documentclass[letterpaper,11pt,runningheads,envcountsame]{llncs}
\usepackage[utf8]{inputenc}
\usepackage[T1]{fontenc}
\usepackage{textcomp}
\usepackage{graphicx}
\usepackage{amssymb,amsfonts,amsmath}
\usepackage[english]{babel}
\usepackage{ae,aecompl} 

\spnewtheorem{conject}[theorem]{Conjecture}{\bfseries}{\itshape}

\begin{document}

\title{More Discriminants\\with the Brezing-Weng Method}
\author{Gaetan Bisson\inst{1,2} \and Takakazu Satoh\inst{3}}
\authorrunning{G. Bisson \and T. Satoh}
\institute{
LORIA, 54506 Vandoeuvre-lès-Nancy, France
\and Technische Universiteit, 5600 Eindhoven, Netherlands
\and Tokyo Institute of Technology, 152-8551 Tokyo, Japan}
\maketitle

\begin{abstract}
The Brezing-Weng method is a general framework to generate
families of pairing-friendly elliptic curves.
Here, we introduce an improvement which can be
used to generate more curves with larger discriminants.
Apart from the number of curves this yields, it provides an easy way
to avoid endomorphism rings with small class number.

\bigskip

\end{abstract}

\section{Introduction}

Since its birth in 2000,
pairing-based cryptography has solved famous open problems in public
key cryptography:
the identity-based key-exchange \cite{sok:cbop},
the one-round tripartite key-exchange \cite{j:aorpftdh}
and the practical identity-based encryption scheme \cite{bf:ibeftwp}.
Pairings are now considered not only as tools for attacking
the discrete logarithm problem in elliptic curves \cite{mov:recliaff}
but as building blocks for cryptographic protocols.

However, for these cryptosystems to be practical, elliptic curves
with an efficiently computable pairing and whose discrete logarithm problem
is intractable are required.

There are essentially two general methods for the generation of such curves:
the Cocks-Pinch method, which generates individual curves,
and the Brezing-Weng method, which generates families of curves
while achieving better $\rho$-values.

Our improvement extends constructions based on these methods
by providing more curves with discriminants larger than
what the constructions would normally provide
(by a factor typically up to $10^9$ given current complexity
of algorithms for computing Hilbert class polynomials).
In the Cocks-Pinch method the discriminant can be freely chosen
so our improvement is of little interest in this case;
however, the Cocks-Pinch method is limited to $\rho\approx 2$.
To achieve smaller $\rho$-values, one has to use the Brezing-Weng method
where known efficient constructions mostly deal with small
(one digit) discriminants; our improvement then provides
an easy and efficient way to generate several curves
with a wide range of discriminants, extending known constructions
while preserving their efficiency (in particular, the $\rho$-value).

The curves we generate, having a larger discriminant,
are possibly be more secure than curves
whose endomorphism ring has small class number
---even though, at the time of this writing,
no attack taking advantage of a small class number is known.
To say the least, our improvement brings a bit of diversity
to families of curves as generated by the Brezing-Weng method.

\bigskip

In Section \ref{sec-2}, we recall the general framework for pairing-friendly
elliptic curve generation. Then, in Section \ref{sec-3}, we present
the Brezing-Weng algorithm and our improvement.
Eventually, in Section \ref{sec-4}, we study practical constructions
and their efficiency; we also present a few examples.

\section{Framework}\label{sec-2}

\subsection{Security Parameters}

Let $\mathcal E$ be an elliptic curve
defined over a prime finite field $\mathbb F_p$.
We consider the discrete logarithm problem in some
subgroup $\mathcal H$ of $\mathcal E$ of large prime order $r$.
In addition, we assume that $r$ is different from $p$.

For security reasons, the size of $r$ should be
large enough to avoid generic discrete logarithm attacks.
For efficiency reasons, it should also not be too small
when compared to the size of the ground field; indeed,
it would be impractical to use the arithmetic of a very large
field to provide the security level that could be achieved
with a much smaller one.
Therefore, the so-called $\rho$-value
\[\rho:=\frac{\log p}{\log r}\]
must be as small as possible.
Note that, for practical applications, it is desirable that the parameters
of a cryptosystem (here, $p$) be of reasonable size relatively to the
security provided by this cryptosystem (here, $r$), which is precisely
what a small $\rho$-value asserts.

\bigskip

We wish to generate such an elliptic curve
and ensure that it has an efficiently computable pairing,
that is a non-degenerate bilinear map from $\mathcal H^2$ to some cyclic group.

Known pairings on elliptic curves, i.e. the Weil and Tate pairings,
map to the multiplicative group of an extension of the ground field.
By linearity, the non-degeneracy of the pairing (on the subgroup
of order $r$) forces the extension to contain primitive
$r^\text{th}$ roots of unity.
Let $\mathbb F_{p^k}$ be the minimal such extension;
the integer $k$ is called the embedding degree.
It can also be defined elementarily as
\[k=\min\{i\in\mathbb N:r\mid p^i-1\}\]

There are different ways of evaluating pairings,
each featuring specific implementation optimizations.
However, all known efficient methods are based on Miller's algorithm
which relies on the arithmetic of $\mathbb F_{p^k}$.
Therefore, the evaluation of a pairing can only be carried out
when $k$ is reasonably small.

\bigskip

In addition, the discrete logarithm problem must be practically intractable
in both the subgroup of the curve and the multiplicative group
of the embedding field. At the time of this writing,
minimal security can be provided by the bounds
\[\log_2 r\geq 160
\text{ and } k\log_2 p\geq 1024\]
However, these are to evolve and, as the bound on $k\log_2 p$
is expected to grow faster than that on $\log_2 r$
(because the complexity of the index-calculus attack on finite fields
is subexponential whereas that of elliptic curve discrete logarithm
algorithms are exponential), we have to consider larger embedding degrees
in order to preserve small $\rho$-values.

\subsection{Curve Generation}

In order to generate an ordinary elliptic curve with a large prime order subgroup
and an efficiently computable pairing,
we look for suitable values of the parameters:
\begin{itemize}
\item $p$, the cardinality of the ground field;
\item $t$, the trace of the Frobenius endomorphism of the curve
(such that the curve has $p+1-t$ rational points);
\item $r$, the order of the subgroup;
\item $k$, its embedding degree.
\end{itemize}

Here, ``suitable'' means that there exists a curve achieving those values.
This consistency of the parameters can be written as the following list
of conditions:
\begin{enumerate}
\item $p$ is prime.
\item $t$ is an integer relatively prime to $p$.
\item $\left|t\right|\leq 2\sqrt{p}$.
\item $r$ is a prime factor of $p+1-t$.
\item $k$ is the smallest integer such that $r\mid p^k-1$.
\end{enumerate}

By a theorem of Deuring \cite{d:dtdmef},
Conditions 1--3 ensure that there exists an ordinary elliptic curve
over $\mathbb F_p$ with trace $t$.
The last conditions then imply that its
subgroup of order $r$ has embedding degree $k$.

When $r$ does not divide $k$
---which is always the case in cryptographic applications
as we want $k$ to be small (for the pairing to be computable)
and $r$ to be large (to avoid generic discrete logarithm attacks)---
Condition 5 is equivalent to $r\mid\Phi_k\left(p\right)$,
which is a much more handy equation;
therefore, assuming Condition 4, it is also equivalent to
\[r\mid \Phi_k\left(t-1\right)\]

\bigskip

To retrieve the Weierstrass equation of a curve with such parameters
using the complex multiplication method,
we need to look at $-D$, the discriminant (which need not be squarefree)
of the quadratic order in which the curve has complex multiplication.
Indeed, the complex multiplication method is only effective when
this order has reasonably small class number.
Due to a result of Heilbronn \cite{h:otcniiqf}, in practice
we ask for $D$ to be a small positive integer.

Writing the Frobenius endomorphism as an element
of the complex multiplication order leads to the very simple
condition
\[\exists y\in\mathbb N, 4p=t^2+Dy^2\]
which ensures that $-D$ is a possible discriminant.
It is referred to as the complex multiplication equation.
Note that, instead of being added to the list, this condition
may supersede Condition 3 as it is, in fact, stronger.

Using the cofactor of $r$, namely the integer $h$
such that $p+1-t=hr$, the complex multiplication equation
can also be written as
\[Dy^2=4p-t^2=4hr-\left(t-2\right)^2\]

Note that if both the above equation considered modulo $r$
and the ``original'' complex multiplication equation hold,
we recover the equation that states
that the curve has a subgroup of order $r$.

\bigskip

Assuming $p>5$, the third condition implies that $p$ divides $t$
if and only if $t=0$; therefore, as $p$ is expected to be large,
we only have to check whether $t\neq 0$.
This condition is omitted from the list below
as it (mostly) always holds in practical constructions;
bear in mind that it is required, though.

Finally, we can summarize the requirements to generate
a pairing-friendly elliptic curve;
we are looking for:
\[\left\{\begin{array}{rl}
p,r & \text{primes} \\
t,y & \text{integers} \\
D,k & \text{positive integers} \\
\end{array}\right.
\text{ such that }
\left\{\begin{array}{l}
r\mid Dy^2+\left(t-2\right)^2 \\
r\mid \Phi_k\left(t-1\right) \\
t^2+Dy^2=4p \\
\end{array}\right.\]

In practical computations, $r$ may not necessarily be given as a prime.
However, if $r$ is a prime times a small cofactor,
replacing it by that prime leads to the generation of a pairing-friendly
elliptic curve without affecting much the $\rho$-value.
Therefore, this slightly weaker condition is acceptable.

\section{Algorithms}\label{sec-3}

Let us fix $D$ and $k$ as small positive integers.
The Cocks-Pinch method consists in solving the above equations
to retrieve values of $p$, $r$, $t$ and $y$; it proceeds in the following way:
\begin{enumerate}
\item Choose a prime $r$ such that the finite field $\mathbb F_r$
contains $\sqrt{-D}$ and $z$, some primitive $k^\text{th}$ root of unity.
\item Put $t=1+z$ and $y=\frac{t-2}{\sqrt{-D}}\mod{r}$.
\item Take lifts of $t$ and $y$ in $\mathbb Z$
and put $p=\frac{1}{4}\left(t^2+Dy^2\right)$.
\end{enumerate}
This algorithm has to be run for different parameters $r$ and $z$ until
the output $p$ is a prime integer; then, the complex multiplication method
can be used to generate an elliptic curve over $\mathbb F_p$ with $p+1-t$ points,
a subgroup of order $r$ and embedding degree $k$.

Asymptotically, pairing-friendly elliptic curves generated
by this algorithm have $\rho$-value $2$.

\subsection{The Brezing-Weng Method}

The Brezing-Weng method starts similarly by fixing small positive integers $D$ and $k$.
Then, it looks for solutions to these equations as polynomials
$p$, $r$, $t$ and $y$ in $\mathbb Q\left[x\right]$.
Once a solution is found, for any integer $x$,
an elliptic curve with parameters
$\left(p\left(x\right),r\left(x\right),
t\left(x\right),y\left(x\right),D,k\right)$
can be generated provided that
$p\left(x\right)$ and $r\left(x\right)$ are prime
and that $t\left(x\right)$ and $y\left(x\right)$ are integers.

To enable this, we expect polynomials $p$ and $r$ to have infinitely many
simultaneous prime values.
There is actually a very precise conjecture on the density of
prime values of a family of polynomials:

\begin{conject}[Bateman and Horn \cite{bh:prbipiov}]
Let $f_1,\dots,f_s$ be $s$ distinct (non-constant) irreducible
integer polynomials in one variable with positive leading coefficient.
The cardinality of $R_N$, the set of positive integers $x$ less that $N$
such that the $f_i\left(x\right)$'s are all prime,
has the following asymptotic behavior:
\[
\operatorname{card} R_N\sim\frac{C\left(f_1,\dots,f_s\right)}{\prod_i\deg f_i}
\int_2^N\frac{du}{\left(\log u\right)^s}
\text{ ~~~~ when }N\rightarrow\infty,
\]
the constant $C\left(f_1,\dots,f_s\right)$ being defined as
\[\prod_{p\in\mathcal P}\left(1-\frac{1}{p}\right)^{-s}\left(1-\frac{1}{p}
\operatorname{card}\left\{x\in\mathbb F_p:\prod_i
f_i\left(x\right)=0\right\}\right)
\]
where $\mathcal P$ denotes the set of prime numbers.
\end{conject}
The latter constant quantifies how much the $f_i$'s differ from
independent random number generators, based on their behavior over finite
fields; it can, of course, be estimated using partial products.

However, if we only need a quick computational way of checking
polynomials $p$ and $r$, we may use a weaker corollary,
earlier conjectured by Schinzel \cite{ss:schclnp}
and known as \emph{hypothesis H}, which just consists in assuming that
the constant $C\left(f_i\right)$ is non-zero.
Consider two polynomials, $p$ and $r$;
in that case, the corollary states that, provided that
\[\gcd\left\{p\left(x\right)r\left(x\right):x\in\mathbb Z\right\}=1\]
the polynomials $p$ and $r$ have infinitely many simultaneous prime values.

Actually, there is a subtle difference with the polynomials we
are dealing with here: they might have rational coefficients.
However, we believe that the hypothesis of the above conjecture
can be slightly weakened as
\[\gcd\left\{p\left(x\right)r\left(x\right):x\in\mathbb Z\text{ such that }
p\left(x\right)\in\mathbb Z\text{ and } r\left(x\right)\in\mathbb Z\right\}=1\]
so to work with families of rational polynomials.
Of course, we use the convention $\gcd\emptyset=0$
(in case there is no $x$ such that
both $p\left(x\right)$ and $r\left(x\right)$ are integers).

\bigskip

Given small positive integers $D$ and $k$, the Brezing-Weng method
works as follows:
\begin{enumerate}
\item Choose a polynomial $r$ with positive leading coefficient
such that $\mathbb Q\left[x\right]/\left(r\right)$
is a field containing $\sqrt{-D}$ and $z$,
some primitive $k^\text{th}$ root of unity.
\item Put $t=1+z$ and $y=\frac{t-2}{\sqrt{-D}}$
(represented as polynomials modulo $r$).
\item Take lifts of $t$ and $y$ in $\mathbb Q\left[x\right]$
and put $p=\frac{1}{4}\left(t^2+Dy^2\right)$.
\end{enumerate}
This algorithm has to be run for different parameters $r$ and $z$ until
the polynomials $p$ and $r$ satisfy the above conjecture.
Then, we might be able to find values of $x$
at which the instantiation of the polynomials yields a
suitable set of parameters and thus generate an elliptic curve.

To heuristically check whether $p$ and $r$ satisfy the above conjecture,
we compute the $\gcd$ of the product $p\left(x\right)r\left(x\right)$
for those $x\in\left\{1,\dots,10^2\right\}$ such that $p\left(x\right)$ and
$r\left(x\right)$ are both integers. If this $\gcd$ is $1$,
the hypothesis of the conjecture is satisfied;
otherwise, we assume it is not.

\bigskip

The main feature of this algorithm is that the $\rho$-value
of the generated curves is asymptotically equal to
$\frac{\deg p}{\deg r}$;
therefore, a good $\rho$-value will be achieved
if the parameters $\left(D,k,r,z\right)$
can be chosen so that the polynomial $p$ is of degree close to that of $r$.
Because of the way $p$ is defined,
the larger the degree of $r$ is, the more unlikely this is to happen.

Such wise choices are rare and mainly concerned with small discriminants;
indeed, when $D$ is a small positive integer, $\sqrt{-D}$
is contained in a cyclotomic extension of small degree
which can therefore be taken as $\mathbb Q\left[x\right]/\left(r\right)$,
thus providing a $r$-polynomial with small degree.

There exist a few wise choices for large $D$
(cf. Paragraph 6.4 of \cite{fst:atopfec})
but those are restricted to a small number of polynomials
$\left(p,r,t,y\right)$ and do not provide as many families
as we would like.

\subsection{Our Improvement}

The key observation is that, if there exists
an elliptic curve with parameters $\left(p,r,t,y,D,k\right)$,
then for every divisor $n$ of $y$
there also exists an elliptic curve with parameters
$\left(p,r,t,\frac{1}{n}y,Dn^2,k\right)$.
Note that this transformation preserves the ground field
and the number of point of the curve, and therefore its $\rho$-value.

For one-shot Cocks-Pinch-like methods, this is of little interest
since we could have set the discriminant to be $-Dn^2$ in the first place.
However, for the Brezing-Weng method where good choices
of the parameters $\left(D,k,r,z\right)$ are not easily found,
it provides a way to generate curves with a wider range of discriminants
with the same machinery that we already have.

\bigskip

This improvement works as follows:
\begin{enumerate}
\item Generate a family $\left(p,r,t,y,D,k\right)$ using the Brezing-Weng method.
\item Choose an integer $x$ such that $p\left(x\right)$ and $r\left(x\right)$ are
prime, and $t\left(x\right)$ and $y\left(x\right)$ are integers.
\item Compute the factorization of $y\left(x\right)$.
\item Choose some divisor $n$ of $y\left(x\right)$ and generate a curve
with parameters $\left(p\left(x\right),r\left(x\right),t\left(x\right),
\frac{1}{n}y\left(x\right),Dn^2,k\right)$ using the complex multiplication method.
\end{enumerate}

In Step 3, we do not actually have to compute the complete factorization
of $y\left(x\right)$. Indeed, $n$ cannot be too large in order for the
complex multiplication method with discriminant $-Dn^2$ to be practical.
So, we only have to deal with the smooth part of $y\left(x\right)$.

However, to avoid efficiently computable isogenies between the original curve
(with $n=1$, as generated by the standard Brezing-Weng method) and our curve,
$n$ must have a sufficiently large prime factor \cite{g:cibecoff}.
Indeed, such an isogeny would reduce the discrete logarithm problem
from our curve to the original curve.

These constraints are best satisfied when $n$ is a prime in some interval.
Specifically, let $D$ be fixed and consider prime values for the variable $n$;
the complexity of computing the Hilbert class polynomial
(with discriminant $-Dn^2$)
is $\operatorname{\Theta}\left(n^2\right)$ \cite{e:tcocpcvfpa}
and that of computing the above-mentioned isogeny
is $\operatorname{\Theta}\left(n^3\right)$ \cite{g:cibecoff}.

Therefore, we recommend to choose a prime factor $n$ of $y\left(x\right)$
as large as possible among those $n$ such that the complex multiplication
method with discriminant $-Dn^2$ is practical, that is,
the Hilbert class polynomial is computable in reasonable time.
Given current computing power, $n\approx 10^5$ seems to be a good choice;
however, to choose the size of the parameter $n$ more carefuly,
we refer to a detailed analysis of the complexity \cite{e:tcocpcvfpa}.

By a theorem of Siegel \cite{s:udcqz}, when $D$ is fixed,
the class number of the quadratic field with discriminant $-Dn^2$
grows essentially linearly in $n$.
Therefore, with $n$ chosen as described above,
the class number is expected to be reasonably large. This helps avoiding
potential (though not yet known) attacks on curves with principal
or nearly-principal endomorphism ring.

\subsubsection{A toy example.}
Let $D=8$, $k=48$ and $r=\Phi_k$ (the cyclotomic polynomial of order $k$).

As $x$ is a primitive $k^\text{th}$ root of unity in
$\mathbb Q\left[x\right]/\left(r\right)$, put
\[t\left(x\right)=1+x\text{ and }\sqrt{-D}=2\left(x^6+x^{18}\right)\]
The Brezing-Weng method outputs polynomials
\[y\left(x\right)=\frac{1}{4}\left(-x^{11}+x^{10}-x^7+x^6+x^3-x^2\right)
\text{ and }p=\frac{1}{4}\left(t^2+Dy^2\right)\]
and the degree of $p$ is such that this family has $\rho$-value $1.375$.

For example, if $x=137$ then
\[p\left(x\right)=12542935105916320505274303565097221442462295713\]
which is a prime number and $r\left(x\right)$ is a prime number as well.
The next step is to factor $y\left(x\right)$ as
\[y\left(x\right)=
-1\cdot 2\cdot 17\cdot 137^2\cdot 229\cdot 9109\cdot 84191\cdot 706631\]
and $n$ can possibly be any product of these factors.

Take for instance $n=17$, which results in discriminant $-2312$
with class number $16$ (as opposed to class number one which
would be provided by the standard Brezing-Weng method, i.e.
with $n=1$). The Weierstrass equation of a curve with parameters
$\left(p\left(x\right),r\left(x\right),t\left(x\right),
\frac{1}{n}y\left(x\right),Dn^2,k\right)$ is given
by the complex multiplication method as
\[\begin{array}{rcl}
Y^2 = X^3
&+& 935824186433623028047894899424144532036848777X \\
&+& 8985839528233295688881465643014243982999429660;
\end{array}\]
this being, of course, an equation over $\mathbb F_{p\left(x\right)}$.

\section{Constructions}\label{sec-4}

We already mentioned that $n$ should have a large prime factor.
To increase chances for $y$ to have such factors,
we seek constructions where $y$ is a nearly-irreducible polynomial,
i.e. of degree close to that of its biggest
(in terms of degree) irreducible factor.

Many constructions based on the Brezing-Weng method can be found
in Section 6 of the survey article \cite{fst:atopfec}.
However, only few involve a nearly-irreducible $y$
(most of those $y$ are divisible by a power of $x$).
Here, we describe a generic construction that is likely to provide
nearly-irreducible $y$'s.

\subsection{Generic Construction}

Fix an odd prime $D$ and a positive integer $k$.

The extension $\mathbb Q\left[x\right]/\left(r\right)$
has to contain primitive $k^\text{th}$ roots of unity;
the simplest choice is therefore to consider a cyclotomic extension.

So, let us put $r=\Phi_{ke}$ for some integer $e$ to be determined.
Let $\zeta_D$ be a primitive $D^\text{th}$ root of unity;
the Gauss sum
\[\sqrt{\left(\frac{-1}{D}\right)D}
=\sum_{i=1}^{D-1}\left(\frac{i}{D}\right)\zeta_D^i\]
shows that, for $\sqrt{-D}$ to be in $\mathbb Q\left[x\right]/\left(r\right)$,
the product $ke$ may be any multiple of $\varepsilon D$ where $\varepsilon=4$
if $-1$ is a square modulo $D$, $\varepsilon=1$ otherwise.

\bigskip

Therefore, we can use the following setting for the Brezing-Weng method:
\begin{enumerate}
\item Choose an odd prime $D$ and a positive integer $k$.
\item Put $\varepsilon=4$ if $-1$ is a square modulo $D$, $\varepsilon=1$ otherwise.
\item Choose a positive integer $e$ such that $\varepsilon D\mid ke$.
\item Choose a positive integer $f$ relatively prime to $k$.
\item Put $r=\Phi_{ke}$, $z=x^{ef}$.
\item Use the expression
\[\sqrt{-D}=x^{\frac{ke}{\varepsilon}}\sum_{i=1}^{D-1}\left(\frac{i}{D}\right)
x^{i\frac{ke}{D}}\mod r\]
for the computation of $y$ in the Brezing-Weng method.
\end{enumerate}

As the latter polynomial is of large degree,
it can be expected to be quite random once reduced modulo $r$.
Therefore, it is likely to be nearly-irreducible
and so the polynomial $y$ given by the Brezing-Weng method
might also be nearly-irreducible.

\bigskip

To support this expectation, we have computed
$\delta:=\frac{\deg m}{deg y}$
where $m$ is the biggest irreducible factor of
$y=\frac{-1}{D}\left(z-1\right)\sqrt{-D}$,
the polynomials for $z$ and $\sqrt{-D}$ being given by the above algorithm.
There are 4670 valid quadruplets
$\left(D,k,e,f\right)\in\left\{1,\dots,20\right\}^4$
(i.e. for which $D$ is an odd prime and $\varepsilon D\mid ke$);
the following table gives the number of valid quadruplets in this range
leading to values of $\delta$ with prescribed first decimal.
\[\begin{array}{|l||c|c|c|c|c|c|c|c|c|c|c|}
\hline
\delta & 0.0 & 0.1 & 0.2 & 0.3 & 0.4 & 0.5 & 0.6 & 0.7 & 0.8 & 0.9 & 1.0 \\
\hline
 & 79 & 27 & 51 & 72 & 26 & 309 & 388 & 320 & 807 & 1127 & 1464 \\
\hline
\end{array}\]
We see that, in this range, more than $70\%$ of valid quadruplets
lead to a $y$-polynomial whose largest irreducible factor
is of degree at least $0.8\deg y$.

\subsection{Examples}\label{sec-4-2}

\subsubsection{Generic Construction.}

Let $D=3$, $k=9$, $e=1$ and $f=4$.

The Brezing-Weng method outputs the polynomials
\[\begin{array}{rcl}
p\left(x\right)&=&\frac{1}{3}\left(x^8+x^7+x^6+x^5+4x^4+x^3+x^2+x+1\right) \\
y\left(x\right)&=&\frac{1}{3}\left(x^4+2x^3+2x+1\right)\\
\end{array}\]
which represent a family of elliptic curves with $\rho$-value $1.333$.

To generate a cryptographically useful curve from this family,
we look for an integer $x$ such that $p\left(x\right)$ is a prime,
$r\left(x\right)$ is nearly-prime and $y\left(x\right)$ is an integer;
we also have to make sure that $p\left(x\right)^k$ and
$r\left(x\right)$ are of appropriate size for both security and efficiency.

Many such $x$'s are easily found by successive trials;
for instance, in the integer interval
$\left[2^{27};2^{28}\right]$, there are $58812$ of them,
which is only $6$ times less than what a pair of independent
random number generators would be expected to achieve
(calculated as $\int_{2^{27}}^{2^{28}}
\log^{-2}$);
for slightly more than a fifth of these, $y\left(x\right)$ 
has a prime factor in the integer interval $\left[10^4;10^6\right]$,
which can therefore be used as $n$ in our algorithm.


For example, let us put $x=134499652$; we have
\[\begin{array}{rcl}
p\left(x\right)&=&35698341005790839038787210375794\backslash\\
                &&985673959363094188344177147207303
\\r\left(x\right)&=&3\cdot 1973357221157926680445163219766947256676055062891
\\y\left(x\right)&=&419\cdot 153733\cdot 1693488567670454571754477
\end{array}\]

If we choose $n=153733$, the discriminant is $-3\cdot 153733^2$
and has class number $51244$; computations give a Weierstrass equation
for the curve:
\[\begin{array}{rcl}
Y^2=X^3&+&18380344310754022726680092877438\backslash \\
&&217394215740605269665898315768997X
\\&+&3541158719057354715243251263604\backslash \\
&&83038157372705450329206494776897
\end{array}\]

\subsubsection{Sporadic Families.}

Our improvement requires families with nearly-irreducible $y$'s
which is why we described a generic construction that is able to generate
such families for various parameters $\left(D,k\right)$.
However, for a few specific parameters, there are sporadic constructions
with good $\rho$-values that also feature nearly-irreducible $y$'s,
and our improvement produces curves with larger discriminants
without changing $\rho$-values.


To illustrate this, let us consider the Barreto-Naehrig family \cite{bn:pfecopo}
which features the optimal $\rho$-value of $1$ for parameters $D=3$, $k=12$ and
\[\begin{array}{rcl}
p\left(x\right) &=& 6^2x^4 + 6^2x^3 + 4\cdot 6x^2 + 6x + 1 \\
r\left(x\right) &=& 6^2x^4 + 6^2x^3 + 3\cdot 6x^2 + 6x + 1 \\
y\left(x\right) &=& 6x^2 + 4x + 1 \\
\end{array}\]


For instance, if $x=549755862066$, we have
\[\begin{array}{rcl}
p\left(x\right)&=&3288379836712499477504831531496220248757101197293
\\r\left(x\right)&=&13\cdot 61\cdot 4146758936585749656374312380967431265034293149
\\y\left(x\right)&=&151579\cdot 11963326366170669619
\end{array}\]

If we choose $n=151579$, the discriminant is $-3\cdot 151579^2$
and has class number $50526$;
computations give a Weierstrass equation for the curve:
\[\begin{array}{rcl}
Y^2=X^3&+&983842331478040932232760138380470085419271212296X
\\&+&2848148112127026939825061113251126889450914939726
\end{array}\]

\section*{Acknowledgements}

The authors would like to thank Pierrick Gaudry for helpful discussions
and Andreas Enge for computing the explicit curve equations found in
Section \ref{sec-4-2}.
Our gratitude also goes to Tanja Lange for her comments and suggestions
on a draft version of this paper.

\bibliography{document}

\end{document}